\documentclass[12pt, reqno]{amsart}
\usepackage{amsmath, amstext, amsbsy, amssymb, amscd}

\newtheorem{theorem}{Theorem}[section]
\newtheorem{lemma}[theorem]{Lemma}

\newtheorem{corollary}[theorem]{Corollary}
\theoremstyle{definition}
\newtheorem{definition}[theorem]{Definition}
\newtheorem{example}[theorem]{Example}

\theoremstyle{remark}
\newtheorem{remark}[theorem]{Remark}

\numberwithin{equation}{section}

\newcommand{\Supp}{{\rm Supp}}

{\vskip-\lastskip\medskip
  \noindent
  {\em #1.}\enspace
  }%
{\qed\par\medskip
  }

\begin{document}
\title[ Algebraic Stein Varieties ] 
              {  Algebraic Stein Varieties } 

\author[Jing  Zhang]{Jing  Zhang}
\address{ Department of Mathematics and Statistics, 
University at  Albany, SUNY, Albany, NY 12222,  USA}
\email{jzhang@albany.edu}

\begin{abstract} 
It is well-known that the associated  analytic  space of 
 an affine variety defined over $\mathbb{C}$
is Stein but the converse is not true, that is, 
an algebraic Stein variety 
is not necessarily affine. 
In  this paper, 
 we  give  sufficient and necessary  conditions 
for an algebraic  Stein variety  to be affine. 
One  of  our results
is that 
an   irreducible  quasi-projective  variety $Y$
defined over $\mathbb{C}$
 with dimension $d$ ($d\geq  1$)
 is affine if and only if 
$Y$ is Stein,  
   $H^i(Y, {\mathcal{O}}_Y )=0$ for all $i>0$
and  $\kappa(D, X)= d$ (i.e.,  
$D$  is  a big  divisor), 
where  $X$ is a projective  variety  containing  $Y$
and $D$ is an effective divisor  with  support   $X-Y$. 
If $Y$ is algebraic Stein but not affine, we also 
discuss  the   possible  
transcendental  degree  of the  nonconstant
regular functions  on   $Y$.  
We  prove  that  $Y$  cannot  have  $d-1$
algebraically   independent  
nonconstant  regular  functions.
The interesting phenomenon
is that the    transcendental   degree  can be  even if the 
dimension of $Y$ is even  and the degree can be odd if 
the 
dimension of $Y$ is odd.     
\end{abstract}

\maketitle


\date{}
\section{Introduction}

We work over complex number field $\mathbb{C}$   and use the terminology
in Hartshorne's book \cite{H1}.

Affine  varieties
(i.e., irreducible  closed  subsets  of  ${\mathbb{C}}^n$  in  Zariski  topology)
  are important in algebraic  geometry.
  Since  J.-P. Serre discovered 
  his  well-known  cohomology  criterion
  (\cite{H2}, Chapter II, Theorem 1.1), 
the  criteria  for  affineness   have  been   investigated 
by     many  algebraic geometers  
(Goodman and Harshorne \cite{GH}; Hartshorne \cite{H2}, Chapter II;  Kleiman \cite{Kl}; Neeman \cite{N}; Zhang\cite{Zh3}).
Corresponding to affine varieties  in algebraic geometry, 
   in  complex  geometry, Stein varieties hold  similar 
 important position. 
 A  complex  space   $Y$
 is Stein if and only if 
 $H^i(Y, G)=0$ for every analytic coherent sheaf $G$ on $Y$
 and all positive integers $i$. Or equivalently, if  $Y$ is  a
 complex analytic     
variety, then
 $Y$ is Stein if and only if  it is both holomorphically 
 convex and holomorphically separable   (\cite{Gu}, Page 143). We say that  $Y$ is  
holomorphically 
 convex if for any discrete sequence $\{y_n\}\subset Y$,
 there is a holomorphic function $f$ on $Y$ such that
 the supremum of the set $\{|f(y_n)|\}$ is 
 $\infty$. $Y$ is  
holomorphically  separable if for every pair $x,y\in Y$,
$x\neq y$,  there is a holomorphic function $f$ on $Y$ such that
$f(x)\neq f(y)$.

It is well-known that the associated  analytic  space of 
 an affine variety defined over $\mathbb{C}$
is Stein 
(\cite{H2}, Chapter  VI,  Proposition  3.1), 
but the converse is not true, that is, 
an algebraic Stein variety 
is not necessarily affine. 
(An algebraic  variety is a quasi-projective  variety,  that  is, 
a Zariski open subset of a projective  variety.)
The  reason is that  there are algebraic varieties which 
have  plenty of global holomorphic functions
but do not have any nonconstant regular functions. 
 J.-P.  Serre  constructed  the first example, 
an algebraic   surface  which is Stein  but not affine
(\cite{H2}, Chapter  VI, Example 3.2).   Based   on   Serre's  
construction  and  the  K$\mbox{\"{u}}$nneth  
formula \cite{SaW},  it is easy to  construct 
higher  dimensional  algebraic  Stein  varieties   which are  
  not  affine.

Then  naturally we would ask: What is the necessary and sufficient condition 
for an  algebraic Stein 
variety  to be  affine? We will answer this question  by showing 
  the following  
Theorem 1.1 and Theorem 1.2.
   In 1988, Neeman  
   proved that a  quasi-affine  
normal variety $U$ of finite type 
over  $\Bbb{C}$  is 
affine if  $U$  is Stein and   the ring 
$\Gamma (U, {\mathcal{O}}_U )$
of
regular  functions on $U$  is a finitely  generated 
$\Bbb{C}$-algebra (\cite{N}, Proposition 5.5).
A quasi-affine variety is a Zariski open subset of an affine 
variety    contained in ${\Bbb{C}}^n$ for some positive 
integer $n$. So  his   theorem is a nice  local
result.
Our  contribution   is
to give  global criteria which work for any algebraic Stein variety $Y$.
We do not assume that $Y$ is contained in  ${\Bbb{C}}^n$.

Let $Y$  be an  irreducible  quasi-projective  variety 
and $X$ be  an  irreducible projective  variety containing $Y$. 
By further blowing up the closed subset  of the boundary $X-Y$, 
we may assume that  $X-Y$ 
is  of  pure codimension  1   and  is
support   of  an  effective  
divisor  $D$  with   simple  normal   crossings. 
Let $\kappa(D, X)$ be   the $D$-dimension (or Iitaka dimension)
(\cite{I2}, Lecture 3) of  $X$, then the Kodaira dimension $\kappa(X)$
is  $\kappa(K_X, X)$, where $K_X$  is the canonical 
divisor of $X$. 
We will use these notations  throughout  this paper. 
In particular, all sheaves are algebraic and  coherent.

By definition,   $\kappa(D, X)$  measures  the  number  of  algebraically 
independent  nonconstant  regular functions  on  $Y$.
Since  $D$  is effective, the dimension of 
linear  space 
 $H^0(X, {\mathcal{O}}_X(mD))$,
denoted by $h^0(X, {\mathcal{O}}_X(mD) )$, is not zero.    Therefore we always have   
$\kappa(D, X)\geq 0$.

\begin{theorem} An  irreducible  algebraic  Stein surface  $Y$   is
affine  if and only if  
$\kappa(D, X)=2$. 
\end{theorem}
Thanks for the referee's comment about Theorem 1.1.
The referee pointed out that  Theorem  1.1 may not be new 
and the proof is standard. We wrote it separately because 
of the following  reasons. First,
we cannot find  any reference  and 
in surface case, both the condition  and proof  
are   different 
from higher dimensional case.  In Theorem 1.1, 
we do not need the condition 
$ H^i(Y, {\mathcal{O}}_Y )=0 $  for all $i>0$ as in Theorem 1.2
for higher dimensional  varieties, 
where  ${\mathcal{O}}_Y$  is the sheaf  of  regular  functions  on  $Y$. Finally, the surface case is  an initial step therefore  is
also   important
 because  our proof for 
higher dimensional varieties is based on mathematical induction.

\begin{theorem} An irreducible algebraic  Stein  variety   $Y$  
of dimension  $d\geq 1$ 
is  affine  if and only if  
$\kappa(D, X)=d$ and  $ H^i(Y, {\mathcal{O}}_Y )=0 $  for all $i>0$. 
\end{theorem}

\begin{theorem} If $Y$ is an irreducible algebraic    Stein  
variety with dimension
 $d\geq 1$,
then 

(1) $\kappa(D, X)\neq  d-1$;

(2) If   $d=2k$, then  $\kappa(D, X)$ can be any even number 
$i$, $i=0,2,4,..., 2k$; 

(3) If  $d=2k+1$, then  
$\kappa(D, X)$ can be  any  odd  number 
$j$,   $j=1, 3, 5,..., 2k-1, 2k+1$. 

\end{theorem}

\begin{definition}An algebraic variety $Y$ is regularly separable if 
for any two distinct points $y_1$ and $y_2$ on 
$Y$, there is a regular
function $f$ on $Y$ such that $f(y_1)\neq f(y_2)$.
\end{definition}

\begin{corollary} If $Y$ is an  irreducible  algebraic   Stein 
 variety with dimension $d\geq 1$,
then  the following conditions  are equivalent

(1)  $Y$  is affine;

(2) $Y$ is regularly separable and   $H^i(Y, {\mathcal{O}}_Y)=0$
for all $i>0$;

(3) $Y$ is regularly separable and  $\Gamma(Y, {\mathcal{O}}_Y)$
(the ring of regular functions on $Y$)
is a finitely generated  $\Bbb{C}$-algebra;

(4) 
$\kappa(D, X)>\max (1, d-2)$ 
  and   $ H^i(Y, {\mathcal{O}}_Y )=0 $  for all $i>0$;

(5)  $\kappa(D, X)= d$,
and $ H^i(Y, {\mathcal{O}}_Y )=0 $  for all $i>0$.

\end{corollary}

This paper is  organized as follows. In Section 2, we will
 prove Theorem 1.1 and some results for surfaces.
  In Section 3, we  will 
prove  Theorem 1.2, Theorem   1.3  and Corollary 1.5.   

\section{  Surfaces}

A fiber space  is a  proper surjective   morphism
$f: V\rightarrow  W$  between  two varieties  $V$ and $W$
such that  the  general   fiber  is  connected.

  $D$
is a big  divisor if 
 $\kappa(D, X)=d$, where 
$d$ is the dimension of $X$.   If $L$  is a  line  bundle  
on a  projective  manifold   $M$, 
it determines   a  Cartier  divisor  $D$.
We  define   $\kappa(L, M)=\kappa(D,  M)$.

\begin{theorem}[{\bf Fujita}] Let 
$M$ and  $S$  be two projective manifolds.  
 Let $\pi: M\rightarrow S$
  be a  fiber space and let $L$ and $H$ be  line bundles  on $M$
  and $S$ respectively. Suppose that  
  $\kappa (H, S)=\dim S$ and that
  $\kappa (aL-b\pi^*(H))\geq 0$
for certain  positive integers $a$, $b$.  Then
$\kappa(L,M)=\kappa (L|_F, F)+ \kappa (H, S)$
for a general  fiber $F$  of  $\pi$.   
\end{theorem}

We   can  freely  adjust  the coefficients of the effective 
divisor $D$  or blow up a closed subvariety in $D$  
(and still denote the divisor as $D$)
without changing
the $D$-dimension because of the following 
 two  properties (\cite{I1};  \cite{Uen}, Chapter II, Section 5).  

Let  $f$: $X'\rightarrow X$
        be a surjective morphism between  two complete
        varieties $X'$ and $X$. Let $D$ be a divisor on 
        $X$ and $E$ an effective divisor on $X'$ such that
        codim$f(E)\geq 2$, then  
        $$\kappa (f^{-1}(D)+E, X')=\kappa (D, X), 
        $$
where  $f^{-1}(D)$   is the reduced  transform of $D$,
defined to be
$f^{-1}(D)=\sum D_i$, $D_i$'s are the   irreducible components
of  $D$. 
       The  $D$-dimension  
       also  does not depend on the coefficients of $D$ 
       if  $D\geq  0$. 
      Precisely,  let $D_1$, $D_2$, $\cdot$$\cdot$$\cdot$, $D_n$
       be  any divisor on $X$ such that for every $i$, $1\leq i \leq n$, 
       $\kappa (D_i, X)\geq 0$, then for   integers 
       $p_1> 0,\cdot\cdot\cdot$, $p_n>0$,  we have  \cite{I3},
       $$\kappa(D_1+\cdot\cdot\cdot+D_n,X)=
       \kappa(p_1D_1+\cdot\cdot\cdot+p_nD_n,X).
       $$
      In particular, if $D_i$'s are irreducible components of
      $D$ and $D$ is effective, then we can change its coefficients
to different  nonzero  positive 
integers
      but   do not change the $D$-dimension.

 The following Lemma  2.2   is a modification of  Goodman's theorem
(\cite{H2}, Chapter 2,  Theorem  4.2). The advantage of  our 
$D$-dimension  version  is  that  we can generalized  it to  higher 
dimensional case \cite{Zh3},
whereas   Goodman's   surface   theorem  (\cite{H2}, Theorem 4.2) 
is not true  for higher dimensional varieties:
If  $Y$  is a  higher dimensional    affine variety    contained  in
a projective  variety   $X$, then the boundary $X-Y$  may not be 
the support  of  any ample  divisor  on  $X$. Goodman proved that only
after further blowing up the  closed  subset on the boundary,
the new boundary  may  be  the support  of an ample  divisor. 
The partial reason is that
the  intersection theory for higher dimensional varieties  
is much more complicated  and  generally an effective  big  divisor has  
no Zariski  decomposition  if  the dimension is higher than 2.

A divisor $D$  of  $X$ is nef  if for every  irreducible 
curve  $C$  on  $X$,  we have  $D\cdot C\geq  0$. 

\begin{lemma}
Let  $Y$  be an irreducible open smooth algebraic surface.
 Let $X$ be a  smooth  projective  surface containing $Y$. Then $Y$ is affine if and only if
the following three conditions hold

(1) Y contains no complete curves;

(2)  The boundary $X-Y$ is connected;

(3) $\kappa (D, X)=2$, where $D$ is an effective divisor 
with support $X-Y$. 
\end{lemma} 
{\it Proof}. We will prove that the boundary $X-Y$
is the support of an ample divisor.

Since $Y$ contains no complete curves,
the  boundary  $X-Y$  cannot  be blown down to a point. 
And $X-Y$  is  of pure codimension 1 since it is connected. 

Write the Zariski decomposition  $D=P+N$, where 
$N$ is effective and negative  definite, $P$  is effective and nef and 
any prime component of  
$N$  does not intersect $P$  \cite{Za}. We may assume that 
both $P$ and $N$ are integral by multiplying a positive 
integer  to the equation since both $P$ and  $N$ are $\Bbb{Q}$ divisors
($D$ is a Weil divisor but $P$ and $N$   may  have  
rational coefficients \cite{Za}). 
Let  $\Supp{D}=\{D_1, D_2, \cdot\cdot\cdot, D_n \}=X-Y$.
Since   $\kappa (D, X)=2$,  $P^2>0$ 
(\cite{Sa1}; \cite {Ba}, Corollary 14.18). 
First we claim that $\Supp{P}=\Supp{D}=X-Y$. If 
$\Supp{P}\neq X-Y$, then there is a prime component, say $D_1$, 
in $X-Y$ such that  $P\cdot D_1>0$ and  $D_1$ is not a component of 
$P$  since $X-Y$ is connected.   Let 
$$Q=mP+D_1,$$
 where $m$ is  a  big  positive integer.
Then  $Q$ is an effective divisor and $\Supp{Q}=\Supp{P}\cup D_1$. 
Since  $P^2>0$,  we may choose $m$ such  that  
$$ Q^2=m^2P^2+2mP\cdot  D_1+ D_1^2>0.  
$$
For every prime  component  $E$ in $P$, since $P$ is nef and 
$D_1$ is not contained in $\Supp {P}$, for  sufficiently large $m$,
 we have 
$$ Q\cdot E=mP\cdot  E+D_1\cdot E\geq 0, \quad \quad 
D_1\cdot  Q=mD_1\cdot P +D_1^2 >0.
$$
Since $Y$  contains no  complete curves, 
any irreducible complete curve outside $X-Y$ intersects 
$X-Y$.  
Thus we get a new effective divisor $Q$ 
such  that  $Q$ is nef and $Q^2>0$. We may replace $P$ by $Q$
and still call it $P$. 
By finitely many such   replacements, we 
can   find an effective nef  divisor $P$     such that 
$P^2>0$ and $\Supp{P}=\Supp{D}=X-Y$.

We  claim  that the boundary $X-Y$ is the support of an ample divisor. 
In fact, the following three conditions imply  the  ampleness: 

(1) $X-Y$ is connected;

(2) $Y$ contains no complete curves;

(3) There is an effective nef divisor $P$ with supp$P=X-Y$ and $P^2>0$.

There is a proof  for the above  claim in \cite{H2}, page 69-71. 
It is also not hard to prove it 
by  computing  intersection numbers
(we may change the coefficients  if necessary).

\begin{flushright}
 Q.E.D. 
\end{flushright}

Notice that the above lemma holds for complete normal surfaces. 
For a complete  normal surface $X$, the intersection theory  is due to
Mumford \cite{Mu2}.
Let Div$(X)$  be the group of  Weil divisors
of $X$.  Let  Div$(X, {\Bbb{Q}})=$Div$(X)\otimes  {\Bbb{Q}}$
be the  group of  $\Bbb{Q}$-divisors. 
The intersection pairing     
$${\mbox{Div}}(X, {\Bbb{Q}})\times {\mbox{Div}}(X, {\Bbb{Q}})
\rightarrow  \Bbb{Q}
$$
is defined  in the following way.  Let $\pi: X'\rightarrow X$
be a resolution and let $A=\cup E_i$
denote the exceptional  set of $\pi$.
For a $\Bbb{Q}$-divisor $D$ on $X$
we define the inverse image 
$\pi^*D$ as
$$ \pi^*D=\bar{D}+\sum a_iE_i 
$$
where $\bar{D}$  is the strict transform  of $D$ by $\pi$
and the rational numbers $a_i$ are uniquely determined by the equations 
$\bar{D}\cdot  E_j+\sum a_iE_i\cdot   E_j=0$  for all $j$. For two divisors 
$D$ and $D'$ on $X$,   their intersection number  is  defined to be   
$$   D\cdot D'=\pi^*D\cdot \pi^*D'.
$$

\begin{lemma}{\bf[Fujita]} Let  $D$ be an effective  $\Bbb{Q}$-divisor
on a normal projective  surface $X$. Then  there exists a unique 
decomposition 
$$D=P+N$$
satisfying the following conditions:

(1) N is an effective $\Bbb{Q}$-divisor and either N=0 or the intersection 
matrix of the irreducible components of N is negative definite;

(2) P is a nef $\Bbb{Q}$-divisor and the intersection of P 
with each irreducible 
component of N is zero. 
\end{lemma}

\begin{lemma}
Let  $Y$  be an irreducible  algebraic  surface.
 Let $X$ be a projective  surface containing $Y$.
 Then $Y$ is affine if and only if
the following three conditions hold

(1) Y contains no complete curves;

(2)  The boundary $X-Y$ is connected;

(3) $\kappa (D, X)=2$, where $D$ is an effective divisor 
with support $X-Y$.   
\end{lemma} 
{\it Proof}. If we have a surjective and finite morphism from 
a variety  $Y'$ to $Y$, then $Y$ is affine if and only if $Y'$ is affine 
by Chevalley's theorem (\cite{H2}, Chapter 2, Corollary 1.5). 
Thus $Y$ is affine if and only if its normalization is affine.
 So we may assume that both 
$Y$ and $X$ are normal by taking their normalization. On a normal 
 projective surface, the intersection theory and Zariski decomposition
remain true by  Lemma 2.3 \cite{Mu2,  Sa2}. Therefore Lemma 2.2 
holds for normal   projective 
surfaces.
In fact, write the  Zariski decomposition
$D=P+N$ by  Lemma 2.3, then 
$P^2>0$ (\cite{Ba}, Corollary  14.18, Page 222).
By the same argument 
as in the proof of  Lemma 2.2, we can find a 
new effective 
nef divisor, still denoted by $P$, such that  
supp$P=X-Y$.   By changing the coefficients of $P$,
we can  find     an ample divisor supported 
in  $X-Y$.   
\begin{flushright}
 Q.E.D. 
\end{flushright}

\begin{remark}  Lemma 2.4 does not hold  for threefolds. 
If $Y$ is a smooth  algebraic  threefold without complete curves
and the boundary $X-Y$ is connected for a smooth completion $X$ of $Y$,
then  we cannot  claim that  $Y$  is affine  by the following two
conditions:

 (1) The 
boundary  $X-Y$  is   of pure codimension 1;

(2)
$\kappa(D, X)=3$.

The reason is that 
in surface case,
the two conditions, i.e., $Y$ contains no complete curves and
$X-Y$ is connected, imply that for any smooth completion $Z$ 
(may be  different from $X$)
of $Y$,
 the boundary $Z-Y$  is  again 
of pure codimension 1. This is of course  not true in higher dimension.
For instance, remove a hyperplane section $H$ and a line
$L$ from ${\Bbb{P}}^3$, where $L$ is not contained in $H$. 
Let $Y={\Bbb{P}}^3-H-L$.  Then $Y$ contains no complete curves.
 Let 
$f:X\rightarrow {\Bbb{P}}^3$
be the blowing up of ${\Bbb{P}}^3$ along $L$. Then $X$ is a  smooth 
projective threefold and  $Y$ is an open subset of $X$. 
Let $D=f^{-1}(H)+E$, where $E$ is the exceptional divisor.
Then by the properties of $D$-dimension,
$\kappa (D, X)=\kappa (H, {\Bbb{P}}^3)=3$. But $Y$ is not affine 
since  the boundary ${\Bbb{P}}^3-Y$ is not of pure codimension
1  (\cite{H2}, Chapter 2, Proposition 3.1).  
 \end{remark}

\begin{corollary} Suppose that we have a surjective morphism from an 
irreducible  smooth   algebraic surface $Y$ to a smooth affine curve $C$.
Let $X$ be a  smooth projective surface containing $Y$.
If      $Y$ contains no complete curves and  
the boundary $X-Y$ is connected, then 
$Y$ is affine.
 \end{corollary}
{\it Proof}. Let $f:Y\rightarrow C$ be the given morphism.
 Then  $f$ gives a rational map from $X$ to $\bar{C}$, where 
$\bar{C}$  is the smooth completion of $C$.  Resolve the
indeterminacy of $f$ on the boundary $X-Y$.
We may replace $X$ by its suitable blowing up and assume that
$f:X\rightarrow  \bar{C}$  is a  surjective and proper morphism. 
Notice that this  procedure  does not  change $Y$.
$Y$ is still an open subset of $X$. 
By Stein factorization, we may assume that every fiber is 
connected and general fiber is smooth.
Pick a point $t_1\in \bar{C}-C$, then 
$$h^1(\bar{C}, {\mathcal{O}}_{\bar{C}}(nt_1))
=0$$
since $nt_1$ is ample for large $n$
(\cite{H1}, Chapter IV, Corollary  3.3). 
By the  Riemann-Roch formula,
$$h^0(\bar{C}, {\mathcal{O}}_{\bar{C}}(nt_1))=
1+n-g(\bar{C}).$$
So 
$\kappa(t_1, \bar{C})=1$ (\cite{Uen}, Chapter II, Section 5).  For a general point 
$t\in C$,
we  may  assume   that  the fiber  $f^{-1}(t)$
determined by  $t$   is smooth   (\cite{Uen}, Corollary 1.8).
 By  Riemann-Roch  formula, there is  a positive integer $m$,
 such that 
$h^0(\bar{C}, {\mathcal{O}}_{\bar{C}}(mt_1-t))>1$. 
 Let $s$ be a nonconstant section of 
$H^0(\bar{C}, {\mathcal{O}}_{\bar{C}}(mt_1-t))$, then
$$ {\mbox{div}} s +mt_1-t\geq 0.
$$
Pull it back to $X$, we have 
$$  f^*({\mbox{div}} s +mt_1-t)={\mbox{div}} f^*(s)+mf^*(t_1)-f^*(t)\geq 0.
$$
Let $D_1=f^*(t_1)$ and $F=f^*(t),$ then 
$f^*(s)$  is a 
nonconstant
section of   $H^0(X, {\mathcal{O}}(mD_1-F))$.
So 
$$h^0(X, {\mathcal{O}}(mD_1-F))>1.$$  
Choose an effective divisor 
$D$ with support $X-Y$ such that  $D_1\leq D$, then we  have
$$ h^0(X, {\mathcal{O}}(mD-F))\geq  h^0(X, {\mathcal{O}}(mD_1-F))>1. 
$$
Since  $F|_Y$ is a smooth affine curve (\cite{H2}, Chapter 2, Proposition 4.1),  and $F$  intersects $D$ on the boundary  $X-Y$,
$D|_F$  is an effective divisor on  $F$.  So 
$$ h^0(F, {\mathcal{O}}_F(mD|_F))\geq  m+1-g(F).
$$
Therefore
 $\kappa(D|_F, F)=1$  (\cite{Uen}, Chapter II, Section 5).
 By Theorem 2.1, 
$$\kappa(D, X)=\kappa(mD, X)=\kappa(mD|_F, F)+\kappa(t, \bar{C})=2.$$
By Lemma 2.2, $Y$ is affine.
\begin{flushright}
 Q.E.D. 
\end{flushright}

It is easy to see that we can drop the smooth condition 
in Corollary 2.6.

\begin{remark} 
Corollary 2.6 does not hold for threefolds. We have the following 
counter-example.

Let $C_t$ be a smooth projective elliptic curve 
defined by $y^2=x(x-1)(x-t)$, $t\neq 0, 1$.  Let $Z$ be the elliptic surface
defined by the same equation. We have surjective morphism from $Z$ to 
$C={\Bbb{C}}-\{0, 1\}$ such that for every $t\in C$, the fiber $f^{-1}(t)=C_t$.   
In \cite{Zh2}, we proved  the following  two claims.

 1. There is a rank 2 vector bundle $E$ on $Z$ such that when 
   restricted to $C_t$, $E|_{C_t}=E_t$ is the unique nonsplit 
   extension of 
    ${\mathcal{O}}_{C_t}$ by ${\mathcal{O}}_{C_t}$, where $f$
     is the morphism from $Z$ to $C$.  

2. There is a divisor $D$ on 
$X={\Bbb{P}}_Z(E)$ such that when restricted to the 
ruled  surface  $X_t={\Bbb{P}}_{C_t}(E_t)$,
$D|_{X_t}=D_t$ is the canonical section of $X_t$.  

  Let  $Y=X-D$, then  $Y$  contains  no complete  curves 
since  we have $H^i(Y, \Omega^j_Y)=0$
for all $i>0$ and $j\geq 0$ \cite{Zh2}.  And  we have the  morphism
from $Y$   to  $C$.  But  $Y$  is not  affine
since  $Y$  has only one algebraically independent 
nonconstant
regular function  \cite{Zh2}.  I do not know  whether
$Y$  is Stein or not.

\end{remark}
\begin{theorem} Let $Y$ be a  Stein algebraic surface
contained in a projective surface $X$. 
Then $Y$ is affine if and only if 
$\kappa(D, X)=2$, where $D$ is an effective divisor with support $X-Y$. 
Moreover, if  $Y$ is not affine,  then
$\kappa(D, X)=0$.
\end{theorem}

{\it Proof}.  We may assume that  $Y$ is normal as 
in the proof of Lemma 2.4.
 Since $Y$ is Stein,  it contains no complete
curves and  the boundary is connected and of pure codimension 1  (\cite {N}, Proposition 3.4).
Now the first claim  is obvious by   Lemma  2.4. 

For the second claim, notice that if $\kappa(D, X)=1$,
then we have surjective morphism from $Y$ to a smooth affine
curve $C$. If  $Y$  is smooth, by Corollary 2.6,  
 $Y$ is affine which contradicts  the fact  
$\kappa(D, X)=1$. 
If  $Y$  is not smooth,  let     $\pi: X'\rightarrow  X$
be  a birational   proper   surjective   morphism   such   that
$X'$  is smooth, then  by  Corollary 2.6,
$\kappa(\pi^*D, X')=2$. Therefore 
$$\kappa(D, X)=\kappa(\pi^*D, X')=2.
$$  
The maximum $D$-dimension   implies that $Y$  is affine  by Lemma 2.4.

Since     $\kappa(D, X)\geq   0$  ($D$  is effective), 
 if  $Y$  is not affine, we have  
$\kappa(D, X)=0$. 

\begin{flushright}
 Q.E.D. 
\end{flushright}

We  have finished  the proof of Theorem 1.1.

\begin{remark}   If  $Y$ is a smooth algebraic 
Stein variety with dimension 3, then $\kappa(D, X)\neq 2$ 
(the following Theorem 3.3) but 
$\kappa(D, X)=1$ is possible \cite{Zh1, Zh2}.  The case $\kappa(D, X)=0$
is a mystery. I do not know the existence of such algebraic variety.
In surface case, J.-P. Serre gave an example (\cite{H2},
Chapter VI, Example 3.2): The open surface
$Y$ is Stein   and  $\kappa(D, X)=0$ (\cite{Ku}, Lemma 1.8).
 \end{remark}

\section{Higher  Dimensional  Varieties}

\begin{theorem}  If $Y$ is an    irreducible  algebraic   
 variety of dimension $d$, then $Y$  is affine  if and  only if  
$Y$ is  Stein,  $\kappa(D, X)=d$ 
and   $ H^i(Y, {\mathcal{O}}_Y )=0 $  for all $i>0$,
where $X$  is  a  completion  of  $Y$,  $D$  is the effective 
  divisor with support   $X-Y$ (the boundary),
  and ${\mathcal{O}}_Y$   is  the sheaf of 
  regular functions on  $Y$.
\end{theorem}

$Proof$.   One direction  is  trivial. 
If $Y$ is affine, then
it is Stein (\cite{H2}, Chapter  VI,  Proposition  3.1)  and
 $H^i(Y,{\mathcal{O}}_Y)=0$   for all  $i>0$
by Serre's affineness criterion.
 The affineness of $Y$ also implies that  $D$ is big since after
 further blowing up the boundary
$X-Y$, it is the support of an ample divisor 
(\cite{H2}, Chapter 2, Section 6, Theorem 6.1).

 We will prove   the  converse:  
$Y$  is affine  if  
$Y$ is  Stein,  $\kappa(D, X)=d$ 
and   $ H^i(Y, {\mathcal{O}}_Y )=0 $  for all $i>0$.

We may assume  that  $Y$  is normal since  $Y$
is affine if and only if its normalization is affine. 
We also may assume $d>2$ since the case $d=1$ is trivial and 
the surface case has been proved in Section 2. 
We will use induction on the dimension of  $Y$.
Assume  that the theorem holds  for all
$(d-1)$-dimensional  irreducible  varieties. 
We  need  to  prove  that  the theorem holds for $d$-dimensional 
 irreducible 
variety $Y$. The idea  of proof is:
with the inductive assumption, for any irreducible  curve on  $Y$, we can find a regular 
function  on  $Y$  such that when restricted to this curve, the function is not a constant. Then we can apply Goodman and 
Hartshorne's quasi-affineness criterion \cite{GH} then 
Neeman's affineness criterion (\cite{N}, Theorem 4.1).

{\bf Step 1}. Any prime  principle  divisor  
$Z=\{y\in Y, f(y)=0, f\in H^0(Y, {\mathcal{O}}_Y)\}$   on
$Y$  satisfies  the three conditions in  Theorem 3.1.

$Proof$. $Z$  is Stein since it is a closed codimension 1 subvariety of 
the Stein variety $Y$.  We have $ H^i(Z, {\mathcal{O}}_Z )=0 $  
for all $i>0$ 
since 
 $ H^i(Y, {\mathcal{O}}_Y )=0 $  for all $i>0$
and  there is a
  short exact sequence
$$ 0\longrightarrow 
 {\mathcal{O}}_Y
\longrightarrow 
 {\mathcal{O}}_Y
\longrightarrow 
 {\mathcal{O}}_{Z}
\longrightarrow 
0,
$$
where the first  map is  defined by  $f$.

Let  $\bar{Z}$  be  the  irreducible closed subvariety
of  codimension 1 on $X$ such that $Z$  is an open subset of  
$\bar{Z}$. Let  $D_0=D|_{\bar{Z}}$  be  the restriction divisor 
on $\bar{Z}$. The boundary  $\bar{Z}-Z$  is  of
pure codimension 1 on $\bar{Z}$ since $Z$  is Stein. 
We need to prove $\kappa(D_0, \bar{Z})=d-1$.

First suppose that  $Y$  and $X$  are smooth.  
The  defining  regular  function $f$ of   $Z$ on $Y$  gives a morphism 
from  $Y$  to $\Bbb{C}$. Let $A=f(Y)$, then $A$  is a smooth affine curve.  
So $f$  defines a rational  map  from   $X$  to ${\Bbb{P}}^1$. 
By Hironaka's  elimination  of  points of   indeterminacy  of  a 
rational  map,
  we have a proper surjective 
morphism  $f'$ from  a smooth projective variety  $X'$  to  ${\Bbb{P}}^1$
such that  $f'|_Y=f$  and  $Y$  is still an open subset of 
$X'$ (i.e., both $X'$ and $X$ are smooth projective varieties
containing $Y$  and $X'$  can be blown down to  $X$).
We still denote the morphism by $f$ and  $X'$  by  $X$.  
So we have the following   
commutative diagram  
\[
  \begin{array}{ccc}
    Y                           &
     {\hookrightarrow} &
    X                                 \\
    \Big\downarrow\vcenter{%
        \rlap{$\scriptstyle{f|_Y}$}}              &  &
    \Big\downarrow\vcenter{%
       \rlap{$\scriptstyle{f}$}}      \\
A       & \hookrightarrow &
{\Bbb{P}}^1,
\end{array}
\]
where  $f$ is proper  and  surjective.

If  the image of  $D$ under $f$   is  a set of finite points, then $Y$
contains  complete  varieties of codimension 1. This is not possible since $Y$  is Stein. Thus  
$f(D)={\Bbb{P}}^1$. 
It is easy to see that  
$\bar{Z}$  is an irreducible  component of the  fiber 
  $f^{-1}(0)$.  By the property of 
  $D$-dimension  in Section 2,
  we may assume  
  $\bar{Z}=f^{-1}(0)$. 
By Corollary  1.8,    \cite{Uen},
a general fiber on  $X$ is smooth.

Let  $X {\stackrel{h}{\longrightarrow}}
\bar{C}   {\stackrel{\alpha}{\longrightarrow}} {\Bbb{P}}^1$
be the Stein factorization, then  
$\bar{C}$
is a smooth projective curve, 
where the first morphism 
$h$  from  $X$  to  $\bar{C}$ is proper and surjective  and  the  second 
morphism $\alpha$ from  $\bar{C}$   to  ${\Bbb{P}}^1$
is a finite map such that every fiber  of  $h$  in  $X$
is connected.  Let  $C=h(Y)$  be  the  image of  $Y$  under the map $h$.
Then  $C$  is  a smooth affine curve.  The above  commutative diagram becomes 
to
\[
  \begin{array}{ccc}
    Y                           &
     {\hookrightarrow} &
    X                                 \\
    \Big\downarrow\vcenter{%
        \rlap{$\scriptstyle{h|_Y}$}}              &  &
    \Big\downarrow\vcenter{%
       \rlap{$\scriptstyle{h}$}}      \\
C        & \hookrightarrow &
\bar{C}.
\end{array}
\]
 
For  any point  $t\in  C$, the  corresponding    
open  fiber  $Y_t=h^{-1}(t)\cap  Y$    is  
Stein since it is a closed
subvariety of the  Stein variety  $Y$. Since   $\bar{Z}=f^{-1}(0)$
is  irreducible,  $h^{-1}(\alpha^{-1}( 0 ))$  is also irreducible and 
 $$h^{-1}(\alpha^{-1}( 0 ))=f^{-1}(0)=\bar{Z},$$ 
where  $\alpha^{-1}( 0 )\in  C$.

Let  $t$  be a general point of  affine curve $C$. 
By Theorem 5.11, \cite{Uen}, 
let $X_t=h^{-1}(t)$,  $D_t=D|_{X_t}$, then
 we have  
$$  d=\kappa(D, X)\leq  \kappa(D_t, X_t) +1\leq  d.
$$
So
$\kappa(D_t, X_t)=d-1$, 
for a general  fiber  $X_t$. 
By upper semi-continuity  theorem, $$\kappa(D_0, \bar{Z})=d-1.$$

If  $Y$  is not smooth, 
we still have the commutative 
diagram 
\[
  \begin{array}{ccc}
    Y                           &
     {\hookrightarrow} &
    X                                 \\
    \Big\downarrow\vcenter{%
        \rlap{$\scriptstyle{h|_Y}$}}              &  &
    \Big\downarrow\vcenter{%
       \rlap{$\scriptstyle{h}$}}      \\
C        & \hookrightarrow &
\bar{C}
\end{array}
\]
such that  $h$  is proper, surjective and 
every fiber of  $h$ in  $X$  is connected.
 Let $\pi: X'\rightarrow X$ be  a proper  surjective  birational 
morphism  such that  $X'$  is  smooth, then we have  a new 
commutative diagram

 \[
  \begin{array}{ccc}
    Y'                           &
     {\hookrightarrow} &
    X'                                 \\
    \Big\downarrow\vcenter{%
        \rlap{$\scriptstyle{g|_Y}$}}              &  &
    \Big\downarrow\vcenter{%
       \rlap{$\scriptstyle{g}$}}      \\
C        & \hookrightarrow &
\bar{C},
\end{array}
\]
where  $g=h \circ \pi$  and  $Y'=\pi^{-1}(Y)$.
Since  
$\pi$ is birational, 
the two   function fields are isomorphic: 
 ${\Bbb{C}}(X)={\Bbb{C}}(X')$.  
Every  fiber  of  $g$ is still  connected
since  every fiber of  $h$  is connected and 
${\Bbb{C}}(X)={\Bbb{C}}(X')$ (\cite{Sh}, Page 139,
the function field 
${\Bbb{C}}(\bar{C})$  is algebraically closed in  ${\Bbb{C}}(X')$).
So  $\kappa(\pi^*D,  X')=\kappa(D, X)=d$  and
$$  d=\kappa(\pi^*D,  X')\leq  \kappa(\pi^*D|_{X'_a}, X'_a)+1\leq  d,
$$
where  $X'_a=g^{-1}(a)$  is the complete  smooth  fiber  in $X'$
(\cite{Uen}, Section 5) for 
a  general  point  $a\in  C$. 
Thus  $\kappa(\pi^*D|_{X'_a}, X'_a)=d-1$,  which means that 
there are  $d-1$  algebraically independent regular functions on 
$X'_a\cap  Y' $. So  $\pi_*X'_a\cap Y=Y_a$ has 
$d-1$  algebraically independent regular functions, that is, 
 $\kappa(D_a,  X_a)=d-1$, where  $X_a=h^{-1}(a)$ and  
$D_a=D|_{X_a}$.  In fact, 
we have  
$$\kappa(D|_{X_a}, X_a)=\kappa(\pi^*(D|_{X_a}), \pi^{-1}(X_a))=
\kappa(\pi^*D|_{X'_a}, X'_a)=d-1.$$

By  upper semi-continuity theorem, we are done.

{\bf Step 2}. If  $H$  is a   hyperplane   such that  $H$ intersects  $X$
with   an irreducible  variety  $\bar{Z}$, then  the restriction of  
$\bar{Z}$  on $Y$:
$Z=\bar{Z}|_ Y$    can  be    defined  by a regular function 
$f$  on  $Y$, that is, $Z=\{y\in Y, f(y)=0, f\in H^0(Y, {\mathcal{O}}_Y)\}$  is  a prime principle divisor on $Y$. 

$Proof$.  Let  $h$  be the irreducible homogeneous   polynomial 
defining  $H$. 
Since  $X$  is normal and projective,  by Bertini's 
theorem for  normal varieties, 
 there is   $H'$,   a different hyperplane defined
by  a  homogeneous  irreducible  polynomial  $h'$
 such  that
$H'$  intersects $X$  with a normal irreducible  subvariety of codimension 
1     \cite{Sei}. Then $h/h'$  is a  rational function on  $X$ and defines 
$H\cap  X=\bar{Z}$.  Since   $\kappa(D, X)=$dim$X$, by the following lemma, 
any rational function
on $X$  can be written as a quotient of two  regular functions 
on $Y$. So   $h/h'=f/g$, where  both $f$   and $g$  are  regular on   $Y$.
Therefore $Z=\{y\in Y, f(y)=0, f\in H^0(Y, {\mathcal{O}}_Y)\}$.

The following lemma is known  and  the  proof can be found in 
\cite{Mo}. 
\begin{lemma}Let $X$ be  normal proper  over 
$\Bbb{C}$.

(1) If there is an $m_0>0$  such that  for all $m>m_0$, 
$h^0(X, {\mathcal{O}}_X(mD))>0$, then 
$$ {\Bbb{C}}(\Phi _{|mD|}(X)) = Q((X, D)),
$$  
where   $\Phi _{|mD|}$  is the  rational 
map  from  $X$  to a projective space 
defined  by   a basis of  $H^0(X, {\mathcal{O}}_X(mD))$.

(2) If  $\kappa(D, X)=$\mbox{dim}X, then  
$\Phi _{|mD|}$ is birational for all $m\gg 0$. 
In particular, ${\Bbb{C}}(X)=Q((X, D))$. 
\end{lemma}
In the above lemma, ${\Bbb{C}}(X)$ is the function field of $X$.
Let $$  R(X, D)=\oplus_{\gamma= 0}^\infty H^0(X, {\mathcal{O}}_X(\gamma D)) 
$$
be the graded $\Bbb{C}$-domain 
and $R^*\subset R$
the  multiplicative  subset  of all nonzero
homogeneous  elements. Then the quotient ring $R^{*-1}R$  
is a graded  $\Bbb{C}$-domain  and its degree 0 part 
$(R^{*-1}R)_0$  is a  field  denoted  by   $Q((X, D))$, i.e., 
$$ Q((X, D))=(R^{*-1}R)_0.
$$

{\bf Step 3}.  For any  irreducible  curve  $F$  on  $Y$, there is a regular 
function    on  $Y$  such that  the restriction  of the function
on  $F$  is not a constant. 

$Proof.$  Let  $\bar{F}$  be the irreducible  complete curve  
on  $X$  containing  $F$  such that   $\bar{F}-F$  is 	 a set of finite points 
on the boundary $X-Y$.

By  theorems  of  Seidenberg \cite{Sei}, 
away from   the finite set  $\bar{F}-F$,
there  is a  hyperplane   $H$  defined by  an irreducible   homogeneous 
polynomial  $h$    such that  $H\cap X$   
is  an irreducible and normal  subvariety  of codimension 1.
Let  $Z=H\cap  Y$, then $Z$  is  irreducible   and normal.
By Step 2, there is a regular function  $f$  on  $Y$  such that
$Z$  is defined by  $f$. We will prove  that
there  is  a regular  function  $r$ on  $Y$
such that 
$r|_F$  is not a constant. 

Since  $H$   is ample and  $H$  does not  contain any point of   $\bar{F}-F$,
the set 
$F\cap  H=F\cap  Z$  is not empty.  If  $F$  is not contained in 
$Z$, then  $f|_F$ is not  a constant.  Suppose that
$F$  is a curve on  $Z$.

By Step 1 and  inductive assumption,
$Z$  is affine.  
From the exact sequence  
$$ 0\longrightarrow 
 {\mathcal{O}}_Y
{\stackrel{f}{\longrightarrow}}
 {\mathcal{O}}_Y
\longrightarrow 
 {\mathcal{O}}_{Z}
\longrightarrow 
0,
$$
we have a surjective map  from 
$H^0(Y, {\mathcal{O}}_Y)$
to  $H^0(Z, {\mathcal{O}}_Z)$  by the vanishing $i$-th cohomology
of  ${\mathcal{O}}_Y$, $i>0$.
Since  $Z$  is affine, there is  a regular function  $r$  on
$Z$  such that $r|_F$  is not a constant. Lift this function to 
$Y$, we are done.

{\bf Step  4}. The  algebraic  Stein variety  $Y$ with 
 $\kappa(D,  X)=d=$dim$Y$
is  quasi-affine, that is, $Y$  is a  Zariski  open subset of 
 an affine variety.

$Proof$. This part  of proof is  due to  Goodman and Hartshorne \cite{GH}.

By a result  of  Goodman  and  Hartshorne \cite{GH}, 
since  the  claim of Step 3 is true, 
there is a proper morphism  $\xi: Y\rightarrow U $
to a quasi-affine  variety  $U$.  Since  
$Y$  is Stein, $Y$  contains no complete  curves. So 
the fiber of the map $\xi$ is of  0 dimensional and
finite. Therefore  $\xi$  is a quasi-finite  morphism.
By Zariski's Main Theorem (\cite{Mu1},  Chapter  III, Section  9),
 $\xi$ factors  through  an open  
immersion  $\alpha: Y\rightarrow  V$
followed by a  finite  morphism  $\beta: V\rightarrow  U$. So 
$Y$  is a quasi-affine  variety.

Since $Y$
 is  quasi-affine and  $H^i(Y, {\mathcal{O}}_Y)=0$  for all $i>0$,
by Neeman's Theorem  (\cite{N}, Theorem 4.1), $Y$  is affine.  

Now we    have  completed  the proof of Theorem 3.1.

\begin{flushright}
 Q.E.D. 
\end{flushright}

\begin{theorem} If $Y$ is an  irreducible quasi-projective  Stein  variety with dimension $d$,
then  
 $\kappa(D, X)\neq  d-1$. 
\end{theorem}
$Proof$.   If  $Y$  is a   curve,
   then 
$Y$ is affine since it is not complete  (\cite{H2}, Chapter II, 
Proposition  4.1;
 \cite{N}).
Choose sufficiently large $n$   such that  $nD$ is ample. The ampleness
implies 
$\kappa(D, X)=\kappa(nD, X)=1$. So the theorem holds for curves.

Suppose that for all  $(d-1)$-dimensional 
varieties the theorem is true.  Let  dim$Y=d>1$.

If $h^0(X, {\mathcal{O}}_X(mD))=1$  for all $m\geq 0$,
then  $\kappa(D, X)=0$. Since $D$ is effective, we
may  assume that   $h^0(X, {\mathcal{O}}_X(mD))>1$
for  all  $m\gg 0$.  Then we have the fiber space 
defined by a nonconstant 
regular function on  $Y$
as
in the proof of Theorem 3.1.  
As  before  we have  the following commutative diagram
\[
  \begin{array}{ccc}
    Y                           &
     {\hookrightarrow} &
    X                                 \\
    \Big\downarrow\vcenter{%
        \rlap{$\scriptstyle{f|_Y}$}}              &  &
    \Big\downarrow\vcenter{%
       \rlap{$\scriptstyle{f}$}}      \\
C        & \hookrightarrow &
\bar{C}
\end{array}
\]
where $C$ is a smooth  affine curve 
embedded in a smooth projective 
curve $\bar{C}$, $f$ is proper and surjective  and
every fiber of $f$ over $\bar{C}$
 is connected. Let  $\pi: X'\rightarrow  X$
be  a proper and surjective   birational  morphism
such that  $X'$ is smooth. Then we have (\cite{Uen}, Chapter II, Theorem 5.13)
$$ \kappa(D, X)= \kappa(\pi^*D, X').
$$
To compute $\kappa(D, X)$, we may assume that  
$X$ is smooth in the above commutative diagram
as  in the proof of   Theorem 3.1, Step 1.

Let $t_1$ be a point in 
$ \bar{C}-C$, then 
$\kappa(t_1, \bar{C})=1$.  For a general point 
$t\in C$, we have
$\kappa(t, \bar{C})=1$. 
 By the Riemann-Roch formula, there is  a positive integer $m$,
 such that 
$$h^0(\bar{C}, {\mathcal{O}}(mt_1-t))>1.$$
  Let $s$ be a nonconstant section of 
$H^0(\bar{C}, {\mathcal{O}}(mt_1-t))$, then
$$ {\mbox{div}} s +mt_1-t\geq 0.
$$
Pull it back to $X$, we have 
$$  f^*({\mbox{div}} s +mt_1-t)={\mbox{div}} f^*(s)+mf^*(t_1)-f^*(t)\geq 0.
$$
Let $D_1=f^*(t_1)$ and $F=f^*(t),$ then  
$$h^0(X, {\mathcal{O}}_X(mD_1-F))>1.$$  
Choose an effective divisor 
$D$ with support $X-Y$ such that  $D_1\leq   D$, then we  have
$$ h^0(X, {\mathcal{O}}_X(mD-F))\geq  h^0(X, {\mathcal{O}}_X(mD_1-F))>1. 
$$
We know that the fibre $F$  is a prime vertical divisor and $f(D)=\bar{C}$.
This implies that  $D|_F$ is an effective divisor
on  $F$.  
Since  $F|_Y$ is a  Stein subvariety of dimension  $d-1$
 (\cite{Gu}, Page 143), by  the  inductive assumption, 
 $\kappa(D|_F, F)\neq  d-2$.
By Fujita's equation, 
$$\kappa(D, X)=\kappa(mD, X)=\kappa(mD|_F, F)+\kappa(t, \bar{C})\neq 
(d-2)+1=
 d-1.$$

\begin{flushright}
 Q.E.D. 
\end{flushright}

\begin{example} Notation is the same as in Theorem 3.3. 

 If   $d=2k$, then  $\kappa(D, X)$ can be any even number 
$i$, $i=0,2,4,..., 2k$.

 $D$ is effective, so  $\kappa(D, X)\geq 0$.
 If $d=2$, then $k=1$  and $\kappa(D, X)$  can be 0 or 2 but not 1
by Corollary  2.8. 
Suppose that the equation holds for $(2k-2)$-dimensional varieties   $V$, $k>1$. 
 We may assume that  $V$  is smooth.
We will construct  an $(2k)$-dimensional variety $Y$  with
the prescribed  $D$-dimension  $i$, an even  number.

 Let  $C$ be an elliptic curve (projective) and $E$ the unique nonsplit 
      extension of $\mathcal{O}$$_C$ by itself.  
      Let ${Z=\Bbb{P}}_C(E)$ and  $F$ be the canonical section,
then $S=Z-F$ is a Stein  surface (\cite{H2}, Chapter VI, Example 3.2)
and  $\kappa(F, Z)=0$  \cite{Ku}. 
Let $A$ be a  smooth affine  curve. 
Let   $Y=V\times S$, 
 $X=\bar{Y}$ be   a smooth completion of $Y$  and $D$ be an
effective divisor  with support  $X-Y$. 
The transcendental degree  of  $H^0(X, {\mathcal{O}}_X(mD))$
does not change since  $S$  has no nonconstant regular functions.
In  fact,  by  K$\mbox{\"{u}}$nneth 
       formula,  

$$\kappa(D, X)=tr.deg_{\Bbb{C}}\oplus_{m\geq 0}
H^0(X, {\mathcal{O}}_X(mD))-1
$$
$$=tr.deg_{\Bbb{C}}\oplus_{m\geq 0}
H^0(\bar{V}, {\mathcal{O}}_{\bar{V}}(mD|_{\bar{V}}))-1
=\kappa(D|_{\bar{V}}, {\bar{V}}),$$  
where $\bar{V}$  is the projective variety  containing  $V$. 
If  $Y=V\times A \times A$, then 
$$\kappa(D, X)=tr.deg_{\Bbb{C}}\oplus_{m\geq 0}
H^0(X, {\mathcal{O}}_X(mD))-1
$$
$$=tr.deg_{\Bbb{C}}\oplus_{m\geq 0}
H^0(X, {\mathcal{O}}_{\bar{V}}(mD|_{\bar{V}}))-1+2
=\kappa(D|_{\bar{V}}, {\bar{V}})+2.$$
By the inductive assumption, the claim holds for $Y$. 
\end{example}

\begin{example}If  $d=2k+1$, then  
$\kappa(D, X)$ can be  any  odd  number 
$j$,   $j=1, 3, 5,..., 2k-1, 2k+1$.

The  calculation is the same as   above. We start with a  Stein 
curve $C$. Let $X=\bar{C}$  be the complete curve 
containing $C$, let $D$ be the boundary  divisor with support  $\bar{C}-C$, then
$\kappa(D, X)=1$. Thus the claim holds for curves.
Suppose that the equation holds for 
$(2k-1)$-dimensional Stein variety $V$, $k\geq 1$. 
Let $Y$ be a Stein variety with dimension $2k+1$.  If 
  we   define    $Y$ 
to be the product space
by adding two
affine  curves, i.e., if  $Y=V\times C\times C$,  then
the dimension of the variety increases by 2 
and the $D$-dimension   also  increases 
  by 2.  
If $Y=V\times S$, then 
the dimension of the variety increases by 2 
but  the $D$-dimension  does not change.
\end{example}

\begin{remark} 1.  I do not know  whether 
there exists an algebraic  Stein  threefold  $Y$  such that  $Y$
has   no  nonconstant   regular   functions, i.e.,
$H^0(Y, {\mathcal{O}}_Y)=\Bbb{C}$. 

2.  If    $Y$ is a quasi-projective  Stein  variety with dimension 
$2k$, then  $\kappa(D, X)\neq  2k-1$.  I do not know whether $\kappa(D, X)=i$,
$i=2k-3, 2k-5,...,  3, 1$  is possible  or  not.

3. Similarly, if $Y$ is a quasi-projective  Stein  variety with dimension 
$2k+1$, then  $\kappa(D, X)\neq  2k$.  I do not know whether $\kappa(D, X)=i$,
$i=2k-2, 2k-4,...,  2, 0$  is possible  or not.

\end{remark}

 Theorem 1.3  has been proved.

\begin{corollary} If $Y$ is an  irreducible  algebraic   Stein 
 variety with dimension $d\geq 1$,
then  the following   conditions   are equivalent

(1)  $Y$  is affine;

(2) $Y$ is regularly separable and   $H^i(Y, {\mathcal{O}}_Y)=0$
for all $i>0$;

(3) $Y$ is regularly separable and  $\Gamma(Y, {\mathcal{O}}_Y)$
is a finitely generated  $\Bbb{C}$-algebra;

(4) 
$\kappa(D, X)>\max (1, d-2)$ 
  and   $ H^i(Y, {\mathcal{O}}_Y )=0 $  for all $i>0$;

(5)  $\kappa(D, X)= d$,
and $ H^i(Y, {\mathcal{O}}_Y )=0 $  for all $i>0$.

\end{corollary}

$Proof$. First, by Theorem  3.3, claim $(4)$  and  $(5)$  are  equivalent.

$(1)\Rightarrow (2)\Rightarrow (5)\Rightarrow (1)$: We only need to prove
$(2)\Rightarrow (5)$, that is, if  $Y$  is regularly separable,  
then  $\kappa(D, X)= d$.

If $Y$  is a curve, then  $Y$ is  affine so 
$\kappa(D, X)= 1$. Suppose the claim holds  for any $(d-1)$-dimensional 
variety. Since  $Y$  is regularly separable, there 
is  a nonconstant regular   
function   $f$  on $Y$. As in the proof of Theorem 3.1,  given by this 
function,
there is 
 a commutative diagram  
\[
  \begin{array}{ccc}
    Y                           &
     {\hookrightarrow} &
    X                                 \\
    \Big\downarrow\vcenter{%
        \rlap{$\scriptstyle{g|_Y}$}}              &  &
    \Big\downarrow\vcenter{%
       \rlap{$\scriptstyle{g}$}}      \\
C        & \hookrightarrow &
\bar{C},
\end{array}
\]
such that  every fiber  $X_t=g^{-1}(t)$  is connected and general fiber is 
irreducible.  As in Step 1, proof of Theorem 3.1, 
we may assume that  $X$  is smooth. Let 
$D_t=D|_{X_t}$.
 By Theorem 2.1, we have 
$$   \kappa(D, X)= \kappa(D_t, X_t)+1=(d-1)+1=d. 
$$

$(3)\Rightarrow (1)$: If $Y$  is  Stein and regularly separable,  then  
  $\kappa(D, X)=d$  so  $Y$  is quasi-affine by  the proof 
of  Step  4, Theorem 3.1. Since 
$\Gamma(Y, {\mathcal{O}}_Y)$
is a finitely generated  $\Bbb{C}$-algebra, $Y$  is affine (\cite{N}, Proposition 5.5).

$(1)\Rightarrow (3)$: Trivial.

The proof is completed.

\begin{flushright}
 Q.E.D. 
\end{flushright}

\end{document}